# ASSESSING EXTREMA OF EMPIRICAL PRINCIPAL COMPONENT FUNCTIONS

By Peter Hall and Céline Vial

*Australian National University*

The difficulties of estimating and representing the distributions of functional data mean that principal component methods play a substantially greater role in functional data analysis than in more conventional finite-dimensional settings. Local maxima and minima in principal component functions are of direct importance; they indicate places in the domain of a random function where influence on the function value tends to be relatively strong but of opposite sign. We explore statistical properties of the relationship between extrema of empirical principal component functions, and their counterparts for the true principal component functions. It is shown that empirical principal component funcions have relatively little trouble capturing conventional extrema, but can experience difficulty distinguishing a "shoulder" in a curve from a small bump. For example, when the true principal component function has a shoulder, the probability that the empirical principal component function has instead a bump is approximately equal to $\frac{1}{2}$. We suggest and describe the performance of bootstrap methods for assessing the strength of extrema. It is shown that the subsample bootstrap is more effective than the standard bootstrap in this regard. A "bootstrap likelihood" is proposed for measuring extremum strength. Exploratory numerical methods are suggested.

**1. Introduction.** The inherent complexity of functional data analysis, as a distinctly infinite-dimensional and infinite-parameter (or nonparametric) problem, means that principal-component methods assume greater importance in FDA than in more traditional, finite-dimensional settings. In particular, there is often no practical opportunity for estimating, in a meaningful way, the "distribution" of a random function. Both the representation of such a distribution, and the slow convergence rates of estimators, throw up obstacles which seem insurmountable in many cases.

Considerations of this type argue that properties of the principal component functions in the distribution of a random function are often going to be of greater importance than properties of the distribution itself. For example, it will be of greater interest to assess peaks and troughs in a principal component function, than to look for extrema in the "density" of the distribution.

The principal component functions appear explicitly in the Karhunen–Loève representation, or expansion, of $X$ (see Section 2.2), where they are weighted by a conventional, scalar random variable. Each part of the function receives the same weight, in terms of the way it contributes to the distribution of $X$. Therefore, places where a principal component function is relatively high, versus places where it is relatively low, have direct interpretation. In particular, local maxima and minima in principal component functions are of explicit importance; they point to places in the domain of a random function where the influence on the function value tends to be relatively strong but of opposite signs. Therefore, identifying extrema in principal component functions, from evidence furnished by their empirical counterparts, is an important part of principal component analysis in FDA.

In some respects this problem is not unlike its counterpart in more classical nonparametric function estimation, where a great deal of effort is often directed toward assessing the numbers of modes and local minima. See, for example, the literature on mode testing and assessment (e.g., [[9, 10, 11, 17, 21, 27, 29, 30, 34], [41, 43]]). For discussion of the mode in nonparametric regression, see, for example, [42, 47, 48]. However, in important respects the two problems are very different. This is reflected in the fact that empirical principal component functions are more accurate estimators of the true principal component functions than conventional nonparametric function estimators are of the true functions. (In particular, they are root-$n$ consistent.) As a result, extrema of empirical principal component functions are more inclined to be close to the correct position than in the case of nonparametric curve estimators.

Moreover, empirical principal component functions are less likely to exhibit spurious "wiggles" in the neighborhood of a real extremum of a true principal component function. This property holds true quite generally, even if the extremum is approached in the manner of a high-degree polynomial. (The extremum of the function with equation $y = x^{2p}$, for a large positive integer $p$, provides an example.) In contrast to these properties, however, empirical principal component functions have considerable difficulty distinguishing a "shoulder" in the true principal component function, from a small "bump" there. (We say that a point $x_0$ is a shoulder point [resp., a bump point] of a continuously differentiable function $f$ if $f'(x_0) = 0$ and $f'(x_0 + x)$ is either strictly positive, or strictly negative [resp., if $f'(x_0) = 0$ and $f'(x_0 + x) f'(x_0 - x) < 0$] for all $x \neq 0$ in a neighborhood of 0. In particular, the origin is a shoulder point of the function $x^{2p+1}$, and a bump point of the function $x^{2p}$.) In such cases the probability that the empirical principal component function also has a shoulder, and the chance that it has instead a single bump, both converge to $\frac{1}{2}$.

In practice, random functions are almost invariably recorded on a discrete grid. (The only exceptions of which we are aware occur in a small number of problems where the functions are recorded by analog means, and even there the data are discretised prior to analysis.) In some contexts the grid is extremely fine; one such example arises in the increasingly common problem of near-infrared spectroscopy, where $X(t)$ denotes the transmission at wavelength $t$, and several thousand values of $t$ are treated at regularly spaced points in an interval. At another extreme, the points at which $X$ is recorded may be rather sparse. For example, in economic data there may be only a dozen values of $X$, representing monthly observations of a process that, in theory at least, operates in the continuum. Particularly in cases such as this the data are smoothed, often by spline methods, prior to obtaining the versions of $X$ to which statistical methodology is applied. Therefore the data may plausibly be supposed to be in the continuum, even if their origin is discrete. In this paper we work with continuous functions, rather than with the discrete information from which those functions are derived.

Bootstrap methods can be used to determine the strength of an extremum, or to assess the possibility that there should be an extremum in a neighborhood of a point where the empirical principal component function has a shoulder. However, it is necessary to employ the subsample bootstrap where the resample size is of smaller order than the sample size (see, e.g., [33]). Using the same size of resample is not as effective. We shall establish all these properties, and suggest a "bootstrap likelihood" for assessing extremum strength. Additionally, we shall develop exploratory numerical methods for addressing the prevalence of extrema and shoulders.

Early development of methodology and theory for principal component analysis of functional data included work of Rao [38], and especially Dauxois, Pousse and Romain [12], who described asymptotic properties of eigenvalues and eigenvectors of sample covariance functions. See also [1, 35, 39, 44, 45]. The technology of FDA has been surveyed and described by Ramsay and Silverman ([37], Chapter 6). There, and more particularly in work of Ramsay and Silverman [36], functional principal component analysis is illustrated by application to real-data examples. Recent work includes that of Cardot [6], Cardot, Ferraty and Sarda [7, 8], Girard [18], James, Hastie and Sugar [26], Boente and Fraiman [4], Huang, Wu and Zhou [23], Mas [28] and He, Müller and Wang [22].

Articles where eigenfunctions play important roles in discrimination and related problems for functional data include those of Huang [24], who used principal components analysis in FDA for gait recognition, and Ferraty and Vieu [15, 16], Glendinning and Herbert [20], Glendinning and Fleet [19] and Biau, Bunea and Wegkamp [3], who employed principal component functions in different ways for classification. Of course, the bootstrap has been used widely in the context of functional data analysis; see, for example, [13, 14, 26, 31, 32, 40, 46].

## 2. Properties of extrema of empirical principal component functions.

2.1. *The case of a fixed distribution of random functions.* First we define functional principal component functions. Let $X_1, X_2, \ldots$ be random functions on a compact interval $\mathcal{I}$, and let $X$ denote a generic $X_i$, with mean $\mu = E(X)$. Put $\bar{X} = n^{-1} \sum_i X_i$,

$$(2.1) \quad K(u,v) = E[\{X(u) - \mu(u)\}\{X(v) - \mu(v)\}] = \sum_{j=1}^{\infty} \theta_j \psi_j(u) \psi_j(v),$$

$$(2.2) \quad \widehat{K}(u,v) = \frac{1}{n} \sum_{i=1}^{n} \{X_i(u) - \bar{X}(u)\}\{X_i(v) - \bar{X}(v)\} = \sum_{j=1}^{\infty} \hat{\theta}_j \widehat{\psi}_j(u) \widehat{\psi}_j(v).$$

The function $K$ is the covariance function of the random process $X$. It can also be interpreted as the kernel of the linear operator that takes a function $\psi$ to $\int_{\mathcal{I}} K(\cdot, u) \psi(v) \, du$. Hence the notation $K$, for kernel, which is common in this setting. The function sequences $\psi_1, \psi_2, \ldots$ and $\widehat{\psi}_1, \widehat{\psi}_2, \ldots$ each comprise an orthonormal basis for the space of square-integrable functions on $\mathcal{I}$. They represent the sequence of "true" principal component functions, and the sequence of empirical principal component functions, respectively. Each pair $(\theta_j, \psi_j)$ in (2.1) represents an (eigenvalue, eigenvector) pair for the linear operator with kernel $K$.

The validity of (2.1) and (2.2) follows from standard results in analysis; see, for example, [25], Chapter 4. The existence of the infinite expansions there is sometimes referred to as Mercer's theorem, although that name is occasionally used for other results. See, for example, [5], Chapter 1. Since $\widehat{K}$, in (2.2), is almost surely a finite rank operator, then the spectral decomposition there is in fact truncated, in the sense that $\hat{\theta}_j$ vanishes for all sufficiently large $j$.

The positive-definiteness property of a covariance function implies that each $\theta_j$ is nonnegative. Therefore, without loss of generality the eigenvalues are ordered so that $\theta_1 \geq \theta_2 \geq \cdots \geq 0$. Likewise, we may assume that $\hat{\theta}_1 \geq \hat{\theta}_2 \geq \cdots \geq 0$. If $\theta_1, \ldots, \theta_{j+1}$ are distinct then the functions $\psi_1, \ldots, \psi_j$ are uniquely defined by (2.1), except that their signs may be reversed. In order to match the sign of $\widehat{\psi}_k$ to that of $\psi_k$ we shall suppose that $\int_{\mathcal{I}} \psi_k \widehat{\psi}_k \geq 0$ for each $1 \leq k \leq j$. Using this convention and assuming that $E(\|X\|^2) < \infty$, it is readily proved that $\|\widehat{\psi}_k - \psi_k\| \to 0$ in probability, and hence that $\int_{\mathcal{I}} \psi_k \widehat{\psi}_k \to 1$ in probability as $n \to \infty$; see, for example, [5], Chapter 4. [Given a continuous, square-integrable function $\psi$ on $\mathcal{I}$, we let $\|\psi\|^2 = \int_{\mathcal{I}} \psi^2$ and $\|\psi\|_\infty = \sup_{u \in \mathcal{I}} |\psi(u)|$.]

The following assumption asks that the $\psi_k$'s have only finite numbers of extrema and horizontal points of inflection (or shoulders), and in particular do not vanish identically on nondegenerate subintervals of $\mathcal{I}$:

(2.3) *For $1 \leq k \leq j$, $\psi_k$ has a Hölder-continuous derivative on $\mathcal{I}$, vanishing at at most a finite number of points $u_{k1} < \cdots < u_{kq_k}$, all of which are interior points of $\mathcal{I}$; and, for $1 \leq k \leq j$ and $1 \leq \ell \leq q_k$, the function $K(u_{k\ell}, \cdot)$ has two square-integrable derivatives on $\mathcal{I}$.*

[If $q_k = 0$ then (2.3) implies that $|\psi_k'|$ is bounded away from zero on $\mathcal{I}_k$.] We insist that each $u_{k\ell}$ be an interior point, since results for points on the boundary have to be framed a little differently. For example, if a local maximum of $\psi_k$ occurs at a point $u_0$ which is an endpoint of $\mathcal{I}$, then the probability that $\widehat{\psi}_k$ has a local maximum in a neighborhood (within $\mathcal{I}$) of that point does not converge to 1. However, this is not the case if $u_0$ is in the interior of $\mathcal{I}$.

The next condition asks that $\psi_k$ behave like a polynomial near its extrema and shoulders, and that it have two smooth derivatives:

(2.4) *For all $1 \leq k \leq j$ and all $1 \leq \ell \leq q_k$, $\psi_k$ has two Hölder-continuous derivatives in a neighborhood of $u_{k\ell}$, and*

$$(\partial/\partial u)^r \{\psi_k(u) - \psi_k(u_{k\ell})\} = (\partial/\partial u)^r A_{k\ell}(u - u_{k\ell})^{r_{k\ell}} + O(|u - u_{k\ell}|^{r_{k\ell} + \eta - r})$$

*as $u \to u_{k\ell}$, where $A_{k\ell} \neq 0$, $r = 0, 1$ or $2$, $r_{k\ell} \geq 2$ is an integer, and $\eta > 0$.*

If $r_{k\ell}$ is even then, in view of (2.4), $u_{k\ell}$ gives a local maximum or local minimum of $\psi_k$ according as $A_{k\ell} < 0$ or $A_{k\ell} > 0$, respectively. If $r_{k\ell}$ is odd then $u_{k\ell}$ is a shoulder-point of $\psi_k$.

A reviewer has reported that in some real-data problems the first eigenfunction is identically constant, or very nearly so. While we have not encountered this ourselves, assumptions such as (2.4) would obviously not be appropriate in such cases.

In Theorem 2.1 below we shall measure the smoothness of the $r$th derivative of $X$ in terms of the finiteness of the moment-based Lipschitz criterion

(2.5) $$\gamma_r(D, \eta) = E\left\{\sup_{u,v \in \mathcal{I}} |X^{(r)}(u) - X^{(r)}(v)|^D |u - v|^{-D\eta}\right\},$$

where $D, \eta > 0$.

Write $\mathcal{I} = [a, b]$, where $-\infty < a < b < \infty$. In the statement of Theorem 2.1 below, let $\varepsilon > 0$ be any positive number not exceeding half the minimum value, over $1 \leq k \leq j$ and $1 \leq \ell \leq q_k + 1$, of $u_{k\ell} - u_{k,\ell-1}$, where $u_{k0} = a$ and $u_{k,q_k+1} = b$. An "$\eta$-neighborhood" of a point $u \in \mathcal{I}$ denotes the set of all real numbers that are distant no more than $\eta$ from $u$. In addition to (2.5) we shall impose conditions which imply that $E|X(u)|^D < \infty$ for each $u \in \mathcal{I}$. This in turn entails $E\|X\|^D < \infty$.

The main aspects of Theorem 2.1 are encapsulated in the following statement: "With probability converging to 1, the extrema of each empirical principal component function, $\widehat{\psi}_k$, correspond exactly to those of the respective true principal component function, $\psi_k$, except that in the neighborhood of a shoulder of $\psi_k$ there may, with probability $\frac{1}{2}$, be exactly two additional extrema. The latter are spurious extrema, in the sense that they arise through stochastic fluctuations and do not reflect actual extrema of $\psi_k$."

THEOREM 2.1. *Assume the eigenvalues $\theta_1, \ldots, \theta_{j+1}$ are distinct, and that with probability 1, $X$ has a Hölder-continuous second derivative on $\mathcal{I}$, with Hölder exponent $\eta > 0$. Suppose too that (2.3) and (2.4) hold. Let $\varepsilon$ be as in the previous paragraph, and suppose $E|X^{(s)}(a)|^D < \infty$ for $s = 0, 1, 2$, and $\gamma_2(D, \eta) < \infty$, for a sufficiently large value of $D > 0$. Then, with probability converging to 1 as $n \to \infty$, and for all $1 \leq k \leq j$ and $1 \leq \ell \leq q_k$, (a)–(c) below hold:* (a) *Within each interval $[u_{k\ell} + \varepsilon, u_{k,\ell+1} - \varepsilon]$, for $1 \leq \ell \leq q_k - 1$, and within each interval $[a, u_{k1} - \varepsilon]$ and $[u_{k,q_k} - \varepsilon, b]$, for $1 \leq k \leq j$, the equation $\widehat{\psi}'_k = 0$ has no solution.* (b) *In each $\varepsilon$-neighborhood of a point $u_{k\ell}$ for which $r_{k\ell}$ is even, $\widehat{\psi}_k$ has exactly one local maximum, or exactly one local minimum, according as $A_{k\ell} < 0$ or $A_{k\ell} > 0$, respectively, and has no horizontal points of inflection.* (c) *In each $\varepsilon$-neighborhood of a point $u_{k\ell}$ for which $r_{k\ell}$ is odd, $\widehat{\psi}_k$ has no more than two extrema. Furthermore:* (d) *When $r_{k\ell}$ is odd, the probability that $\widehat{\psi}_k$ has no extremum, and the probability that there are exactly two extrema (a local maximum or a local minimum, resp.) in an $\varepsilon$-neighborhood of $u_{k\ell}$, both converge to $\frac{1}{2}$.* (e) *If $r_{k\ell}$ is even, then the extremum $\hat{u}_{k\ell}$, say, of $\widehat{\psi}_k$ in a neighborhood of $u_{k\ell}$ [see (b)], satisfies*

$$n^{1/\{2(r_{k\ell}-1)\}}(\hat{u}_{k\ell} - u_{k\ell}) \to N^{1/(r_{k\ell}-1)}$$

*in distribution, where $N$ has a normal $\mathrm{N}(0, \sigma_{k\ell}^2)$ distribution and $\sigma_{k\ell} > 0$.* (f) *If $r_{k\ell}$ is odd, then conditional on $\widehat{\psi}_k$ having two extrema, $\hat{u}_{k\ell}^+ > \hat{u}_{k\ell}^-$ say, in a neighborhood of $u_{k\ell}$ [see (d)], they satisfy*

$$n^{1/\{2(r_{k\ell}-1)\}}(\hat{u}_{k\ell}^+ - u_{k\ell}, \hat{u}_{k\ell}^+ - u_{k\ell}) \to (|N|^{1/(r_{k\ell}-1)}, -|N|^{1/(r_{k\ell}-1)})$$

*in distribution, where $N$ is as in* (e).

2.2. *The case of locally perturbed distributions of random functions.* Results similar to Theorem 2.1 can be obtained in the case of slight alterations to a population, thereby giving insight into the way in which an empirical principal component function $\widehat{\psi}_j$ responds to small bumps in the true principal component function $\psi_j$. For the sake of brevity we shall only summarize the results below, rather than state them as formal theorems. We shall refer to the small alterations as "local perturbations," where the word "local" is used in the sense of the standard terminology "local hypothesis," not in the sense in which "local" is interpreted in statistical smoothing.

In general, the effect of adding a perturbation does not alter the results described in Section 2.1, if the perturbation is of smaller order than the size, $n^{-1/4}$, of noise. On the other hand, the effect of noise on the perturbation is negligible if $n^{-1/4}$ is of smaller order than the size of the perturbation. We shall discuss these properties in detail in a particular case, where the added perturbation is one of the orthogonal functions themselves. There the theory is particularly simple and transparent.

Let $\psi_j(u) = A(u-u_0)^3$, where $u_0$ denotes the midpoint of $\mathcal{I}$, and $A \neq 0$. Thus, $u_0$ is the site of a shoulder of $\psi_j$. We shall add a small bump, of vertical height $\delta$, at $u_0$. Let $\psi_p$, for $p \neq j$, be symmetric about $u_0$, have two continuous derivatives on $\mathcal{I}$, and satisfy $\psi_p \geq 0$ and $\psi_p'' < 0$ on $\mathcal{I}$, together with $\psi_p'(u_0) = 0$. The restrictions on $\psi_j$ and $\psi_p$ ensure that these functions are orthogonal; they are of course easily rendered orthonormal by rescaling. (Note that we are adding a bump function, $\psi_p$, to a shoulder function, $\psi_j$, and in particular are not adding a shoulder function to a shoulder function.)

Note particularly that under the above conditions, and for all sufficiently small $|\delta| > 0$, the function

$$(2.6) \qquad \phi_j \equiv (1-\delta^2)^{1/2}\psi_j + \delta\psi_p$$

has exactly two extrema, one of them at $u_0$ and the other distant $O(|\delta|)$ from $u_0$, giving either a local maximum or a local minimum. We shall explore the extent to which the resulting bump in the local perturbation, $\phi_j$, of $\psi_j$ is visible in the empirical principal component function $\widehat{\phi}_j$.

A random function $X$ whose covariance admits the spectral expansion (2.1) also has the expansion

$$(2.7) \qquad X(u) - E\{X(u)\} = \sum_{k=1}^{\infty} \xi_k \psi_k(u),$$

where the random variables $\xi_k = \int(X - EX)\psi_k$ are uncorrelated and have zero mean, with $E(\xi_k^2) = \theta_k$. Result (2.7) is sometimes referred to as the Karhunen–Loève expansion of $X - E(X)$, although that name is also used for results such as the expansion at (2.1).

We shall take $E(X) = 0$ for simplicity, and perturb $X$ to $Y$, where

$$Y(u) = \sum_{k=1}^{\infty} \xi_k \phi_k(u),$$

with $\phi_k = \psi_k$ for $k \notin \{j, p\}$, $\phi_j$ as at (2.6), $\phi_p = \delta \psi_j - (1-\delta^2)^{1/2} \psi_p$, $0 < |\delta| < 1$ and $\xi_k$ exactly as at (2.7). The corresponding spectral decomposition of the covariance of $Y$ is

$$L(u,v) = \text{cov}\{Y(u)Y(v)\} = \sum_{k=1}^{\infty} \theta_k \phi_k(u) \phi_k(v)$$

$$= K(u,v) + \delta(\theta_j - \theta_p)\{\psi_j(u)\psi_p(v) + \psi_p(u)\psi_j(v)\} + O(\delta^2)$$

as $\delta \to 0$, where $K(u,v)$ and $\theta_1, \theta_2, \ldots$ are as at (2.1).

Let the perturbed random functions $Y_1, \ldots, Y_n$ be independent and identically distributed as $Y$, and let $\delta = \delta(n)$ depend on $n$, converging to zero as $n \to \infty$. Taking $\bar{Y} = n^{-1} \sum_i Y_i$ and

(2.8) $$\widehat{L}(u,v) = \frac{1}{n} \sum_{i=1}^{n} \{Y_i(u) - \bar{Y}(u)\}\{Y_i(v) - \bar{Y}(v)\},$$

being the analogue of $\widehat{K}$ defined at (2.2) but for the random functions $Y_i$ rather than $X_i$, we may obtain a spectral decomposition in the usual way:

(2.9) $$\widehat{L}(u,v) = \sum_{k=1}^{\infty} \tilde{\theta}_k \tilde{\phi}_k(u) \tilde{\phi}_k(v),$$

where $\tilde{\theta}_k$ and $\tilde{\phi}_k$ may be interpreted as estimators of $\theta_k$ and $\phi_k$, respectively. The version of Theorem 2.1 pertaining to the data $Y_i$, rather than $X_i$, includes, among other things, the following properties.

THEOREM 2.2. *Assume the conditions of Theorem 2.1. If $\varepsilon > 0$ is not greater than half the length of $\mathcal{I}$, and $\delta = \delta(n)$ satisfies $n^{1/4}\delta \to \infty$, then with probability tending to 1 as $n \to \infty$, $\tilde{\phi}_j$ has just two extrema, and these occur in an $\varepsilon$-neighborhood of $u_0$; and if $\delta(n) = o(n^{-1/4})$, then the probability that $\tilde{\phi}_j$ has just two extrema in the $\varepsilon$-neighborhood, and the chance that it has no extrema there, both converge to $\frac{1}{2}$. When $\delta(n) = \text{const}.n^{-1/4}$, the probability that there are just two extrema in the $\varepsilon$-neighborhood converges to a number strictly between 0 and 1. For any of these choices of $\delta$, the probability that $\tilde{\phi}_j$ has a horizontal point of inflection in $\mathcal{I}$ equals 0.*

Recall that $\phi_j$ has exactly two extrema, for all sufficiently small values of $|\delta|$. Theorem 2.2 implies that, if $|\delta|$ is of strictly larger size than $n^{-1/4}$, then with probability tending to 1 as $n \to \infty$, the empirical principal component function $\widehat{\phi}_j$ has exactly this number of extrema, and so correctly indicates the presence of the bump in $\phi_j$. However, if $|\delta|$ is of smaller size than $n^{-1/4}$ then the probability that the added bump is reflected in $\widehat{\phi}_j$ converges only to $\frac{1}{2}$. Hence, in this case the bump is not so readily visible. In this sense the "bump" in $\phi_j$ can be detected reliably if it is of larger order than $n^{-1/4}$, but not otherwise.

Reflecting these results, it is not difficult to show that no method is capable of reliably distinguishing between the case where no bump is present (i.e., $\delta = 0$), and that where a bump of size $cn^{-1/4}$, for arbitrarily small but positive $|c|$, has been added. To better appreciate this point, let us assume that the process $X$ is Gaussian, or equivalently, that the random variables $\xi_k$ in (2.7) are independent and normally distributed, with zero means and respective variances $\theta_k$ satisfying $\sum_k \theta_k < \infty$. Fix a sequence $\theta_1 > \theta_2 > \cdots > 0$ with this property, and fix also the eigenvectors $\psi_1, \psi_2, \ldots$, choosing $\psi_j$ and $\psi_p$ as above. More generally, adopt the earlier construction of local perturbations, where a single parameter $\delta$ is involved. Consider a set of just two distributions $D$ of $Y$, the first, $D_0$ say, corresponding to $\delta = 0$, and the second, $D_c$, depending on $c \neq 0$, corresponding to $\delta = cn^{-1/4}$. For a given value of $\alpha \in (0,1)$, let $\mathcal{T}_\alpha$ denote the class of all decision rules $T = T(n)$ that are measurable functions of the data $\mathcal{Y} = \{Y_1, \ldots, Y_n\}$, which can be used to classify the distribution of $Y$ as either $D_0$ or $D_c$, and which satisfy

$$P_{D_0}(T \text{ classifies } D \text{ as } D_0) \geq \alpha.$$

Then,

(2.10) $$\limsup_{c \to 0} \limsup_{n \to \infty} \sup_{T \in \mathcal{T}_\alpha} P_{D_c}(T \text{ classifies } D \text{ as } D_c) \leq 1 - \alpha.$$

This result may be paraphrased by saying that "no test that is capable of accurately detecting $D_0$, can also accurately detect $D_c$ for arbitrarily small values of $|c|$."

**3. Bootstrap-based assessment of extrema of empirical principal component functions.** We begin by describing the $m$-out-of-$n$ bootstrap in the context of estimating principal component functions. Let $m \leq n$, and conditional on $\mathcal{X} = \{X_1, \ldots, X_n\}$, draw a bootstrap resample, $\mathcal{X}^* = \{X_1^*, \ldots, X_m^*\}$, by sampling randomly, with replacement, from $\mathcal{X}$. Define $\widehat{K}^*$ analogously to $\widehat{K}$, at (2.2), and construct its spectral expansion,

$$\widehat{K}^*(u,v) = \frac{1}{m} \sum_{i=1}^{m} \{X_i^*(u) - \bar{X}^*(u)\}\{X_i^*(v) - \bar{X}^*(v)\}$$

$$= \sum_{j=1}^{\infty} \hat{\theta}_j^* \widehat{\psi}_j^*(u) \widehat{\psi}_j^*(v).$$

Of course, $\widehat{\psi}_j^*$ is the bootstrap version of $\widehat{\psi}_j$. We shall use the stochastic variability of extrema of $\widehat{\psi}_j^*$ to measure the relative "strengths" of extrema of $\psi_j$.

The key to our suggested method is the following property. If $m$ is chosen large, but small relative to $n$ [i.e., in asymptotic terms, if $m = m(n) \to \infty$ but $m/n \to 0$], then the numbers of modes of $\widehat{\psi}_k^*$, for $1 \leq k \leq j$, accurately reflect those of $\widehat{\psi}_k$ when the latter are viewed in an unconditional sense. That is, with high probability a "true" extremum in $\psi_k$ will produce exactly one extremum in $\widehat{\psi}_k^*$; but a shoulder point in $\psi_k$ will (in asymptotic terms) produce either no extrema in $\widehat{\psi}_k^*$, or exactly two extrema there, each of these outcomes occurring approximately 50% of the time. Therefore, the relative frequencies with which bumps appear at different places in $\widehat{\psi}_k^*$ give a good guide to the "likelihoods" that real bumps are present in $\psi_k$.

One can construct an informal definition of the likelihood attached to an extremum being in a particular region, as follows. Taking the estimator of $\widehat{\psi}_k$ as a guide, first determine a subinterval, $\mathcal{J}$ say, of $\mathcal{I}$ where a single extremum might lie; replicate the function $\widehat{\psi}_k^*$ and count the proportion, $\hat{p}$, of times (for a given dataset $\mathcal{X}$) that $\widehat{\psi}_k^*$ has at least one extremum in $\mathcal{J}$; and take

$$\hat{\pi} = \max\{0, 2(\hat{p} - \tfrac{1}{2})\} \tag{3.1}$$

to be a measure of the likelihood that there is an extremum of $\psi_k$ in $\mathcal{J}$. Theorem 3.1 asserts that if there is just one point in the interior of $\mathcal{J}$, and no point on the boundary of $\mathcal{J}$ where $\psi_k'$ vanishes, then $\hat{\pi} \to 0$ or $1$ according as that point is a shoulder point, or a proper extremum, respectively.

However, the results are rather different if $m = n$, or more generally if $m/n$ does not converge to zero. There, while it remains true that (with probability close to 1) a "true" extremum in $\psi_k$ will produce exactly one extremum in $\widehat{\psi}_k^*$, it does not hold that a shoulder point in $\psi_k$ will produce extrema in $\widehat{\psi}_k^*$ approximately half the time. Indeed, if $m/n$ does not converge to zero then the proportion of times, $\hat{p}$, that (in the vicinity of a shoulder point of $\psi_k$) there is at least one extremum of $\widehat{\psi}_k^*$, does not converge in probability to a limit. In this sense, the standard, $n$-out-of-$n$ bootstrap is not consistent. Nevertheless, $\hat{p}$ does converge in distribution, to a random variable supported on $[0, 1]$ and distributed symmetrically about $\tfrac{1}{2}$. Intuitively, the reason for this property is that if $m/n$ does not converge to zero then the difference between $\widehat{\psi}_k^*$ and $\widehat{\psi}_k$, which is of size $m^{-1/2}$, is not of strictly larger order than the difference between $\widehat{\psi}_k$ and $\psi_k$, which is of size $n^{-1/2}$. As a result, the latter difference plays a significant role in determining asymptotic properties of the conditional distribution of extrema of $\widehat{\psi}_k^*$.

In the theorem below, given an interval $\mathcal{J} \subseteq \mathcal{I}$, and keeping the sample $\mathcal{X}$ fixed, let $\hat{p}_k^{(\nu)}(\mathcal{J})$ denote the proportion of the resamples $\mathcal{X}^*$ for which there are exactly $\nu$ solutions, in $\mathcal{J}$, of the equation $(\widehat{\psi}_k^*)' = 0$. Let $\varepsilon > 0$ be as defined in the paragraph immediately prior to Theorem 2.1, and recall that $A_{k\ell}$ denotes the constant appearing in (2.4).

THEOREM 3.1. *Assume the conditions of Theorem* 2.1, *and also that* $n \geq m = m(n) \geq n^c$, *for some* $c \in (0, 1)$, *and that* $1 \leq k \leq j$. *Assume initially that* $m/n \to 0$. *Then:* (a) *If* $\mathcal{J}$ *denotes either* $[u_{k\ell} + \varepsilon, u_{k,\ell+1} - \varepsilon]$, *for* $1 \leq k \leq j$ *and* $1 \leq \ell \leq q_k - 1$, *or* $[a, u_{k1} - \varepsilon]$ *or* $[u_{k,q_k} - \varepsilon, b]$, *for* $1 \leq k \leq j$, *then* $\hat{p}_k^{(0)}(\mathcal{J}) \to 1$ *in probability as* $n \to \infty$. (b) *If* $\mathcal{J}$ *denotes an* $\varepsilon$-*neighborhood of a point* $u_{k\ell}$ *for which* $r_{k\ell}$ *is even, then* $\hat{p}_k^{(1)}(\mathcal{J}) \to 1$ *in probability, and in fact the proportion of the solutions of* $(\widehat{\psi}_k^*)' = 0$ *that give a local maximum converges to* 1, *or to* 0, *according as* $A_{k\ell} < 0$ *or* $A_{k\ell} > 0$, *respectively.* (c) *If* $\mathcal{J}$ *denotes an* $\varepsilon$-*neighborhood of a point* $u_{k\ell}$ *for which* $r_{k\ell}$ *is odd, then* $\hat{p}_k^{(\nu)}(\mathcal{J}) \to \frac{1}{2}$ *in probability for* $\nu = 0, 2$. *Next assume that* $m/n \to \rho$, *where* $0 < \rho \leq 1$. *Then* (a) *and* (b) *above continue to hold, but* (c) *should be changed to:* (c)' *If* $\mathcal{J}$ *denotes an* $\varepsilon$-*neighborhood of a point* $u_{k\ell}$ *for which* $r_{k\ell}$ *is odd, then*

$$(\hat{p}_k^{(0)}(\mathcal{J}), \hat{p}_k^{(2)}(\mathcal{J})) \to (\Phi(\rho^{1/2} N), 1 - \Phi(\rho^{1/2} N))$$

*in distribution, where* $\Phi$ *denotes the standard normal distribution function and the random variable* $N$ *has the* $N(0, 1)$ *distribution.*

Of course, if $\rho = 1$ then $\Phi(\rho^{1/2} N)$ has the uniform distribution on $[0, 1]$.

Set up the local perturbation problem as in Section 2.2. In particular, $\psi_j$ denotes a function with a shoulder at $u_0$ (the midpoint of $\mathcal{I}$), and $\psi_p$ is a concave bump, of which a small multiple, $\delta = \delta(n)$, is added to $\psi_j$ to form the locally perturbed version, $\phi_j$, of $\psi_j$; see (2.6). As in Section 2, the principal component function $\phi_j$ is estimated by $\tilde{\phi}_j$, obtained from the spectral decomposition (2.9) of $\widehat{L}$, the latter defined at (2.8) in terms of the dataset $\mathcal{Y} = \{Y_1, \ldots, Y_n\}$. Let $\widehat{L}^*$ denote the bootstrap version of $\widehat{L}$, obtained by computing $\widehat{L}$ not from the original dataset $\mathcal{Y}$ but from a bootstrap resample $\mathcal{Y}^* = \{Y_1^*, \ldots, Y_m^*\}$ drawn by sampling randomly with replacement from $\mathcal{Y}$. Assume that $m = m(n) \geq n^c$ for some $c \in (0, 1)$, and that $m/n \to 0$ as $n \to \infty$. Let $\mathcal{J}$ denote an $\varepsilon$-neighborhood of $u_0$, where $\varepsilon$ is any positive number not greater than half the length of $\mathcal{I}$; and define $\hat{p}_j^{(\nu)}(\mathcal{J})$ to equal the proportion of times, conditional on $\mathcal{Y}$, that the equation $(\widehat{\phi}_j^*)' = 0$ has exactly $\nu$ solutions in the interval $\mathcal{J}$. Assume the conditions of Theorem 2.1. Then the following analogue of a portion of Theorem 2.2 holds:

(3.2) *If* $\delta = \delta(n)$ *satisfies* $m^{1/4} \delta \to \infty$ *then* $\hat{p}_j^{(2)}(\mathcal{J}) \to 1$ *in probability, whereas if* $\delta(n) = o(m^{-1/4})$ *then* $\hat{p}_j^{(\nu)}(\mathcal{J}) \to \frac{1}{2}$ *in probability for* $\nu = 0, 2$. *For either choice of* $\delta$, *the probability that* $\widehat{\phi}_j^*$ *has a horizontal point of inflection in* $\mathcal{I}$ *equals* 0.

Expressed another way, (3.2) states that if $m^{1/4}\delta \to \infty$ then the "bootstrap likelihood" $\hat{\pi}$, at (3.1), will converge to 1 when applied to the interval $\mathcal{J}$ and to the $j$th principal component function. That is, the bootstrap likelihood will correctly signal that a bump, in the form of the function $\psi_p$, has been added to the shoulder at $u_0$ in $\psi_j$. On the other hand, if $\delta(n) = o(m^{-1/4})$ then the bootstrap likelihood will converge to zero, indicating that the added bump has been missed. Therefore, the bump may not be detected if $m$ is too small.

In some respects it might be satisfying to have a purely empirical rule for choosing $m$. However, the results above argue that this is neither practical nor, in the main, actually desirable. We know from (2.10) that no empirical rule can distinguish the case where there is no added bump, from that where a bump of size $n^{-1/4}$ is added; and the results above show that a bootstrap-based test can, in asymptotic terms, distinguish a bump of any order, $\delta(n)$ say, that is strictly greater than $n^{-1/4}$. Indeed, we should choose $m$ so that $m/n \to 0$ but $\delta^4 m \to \infty$. However, in order to achieve this level of sensitivity in a purely empirical way, without an external source of information about bump size, we need to do in advance essentially that which are trying to do now—we need to use empirical evidence to approximate $\delta$ so we can choose $m$ in order to determine empirically whether a bump, of size $\delta$, is present.

The circularity of this argument, and the fact that [in view of (2.10)] it is virtually impossible to accurately estimate $\delta$ when the bumps are small but barely detectable, means that empirical choice of $m$ is not a practical option. Instead, the problem should be addressed from an exploratory angle, for example, starting with $m = n$ and gradually decreasing this quantity. We expect that as $m$ decreases the stochastic fluctuations inherent in the bootstrap will play an increasing role, generally giving rise to more extrema in the functions $\hat{\psi}_k^*$. In consequence, if an extremum is not genuine, $\hat{\pi}$ will tend to increase with decreasing $m$.

## 4. Numerical properties.

4.1. *Application to real data on gait cycle.* This example concerns child gait data, studied by Olshen, Biden, Wyatt and Sutherland [31] and Rice and Silverman [39]. These data consist of records of the angle of the knee and the hip during a gait cycle, recorded at 20 equally spaced time points, for 39 children aged approximately five years. We shall focus on the hip data.

Figure 1 illustrates the first and second empirical principal component functions for one cycle. (To provide a good view of the extrema we represent the data in the interval $[-4, 15]$.) These two eigenfunctions correspond to 82% of total variability, and can be interpreted as follows: The first eigenfunction represents an overall shift with respect to the mean curve, and the second corresponds to an increased angle for the first part of the cycle (in the interval $[0, 8]$), and a delay in the second part.

TABLE 1
*Values of $\hat{\pi}$ for each extremum of the first two principal component functions (PCFs) for the hip-movement component of the child gait data, obtained using 500 bootstrap iterations and different bootstrap sample sizes $m$. The interval $\mathcal{J}$ for each extremum $u_{kl}$ was $[u_{kl} - \varepsilon, u_{kl} + \varepsilon]$*

| | 1st empirical PCF | | | | 2nd empirical PCF | | | |
|---|---|---|---|---|---|---|---|---|
| Extrema | $-1.5$ | 0.5 | 2.5 | 10.5 | 0.5 | 2.5 | 4.5 | 12.5 |
| $\varepsilon$ | 1 | 1 | 1 | 4 | 1 | 1 | 1 | 4 |
| $m=10$ | 0.764 | 0.144 | 0.008 | 0.976 | 0.476 | 0 | 0 | 0.82 |
| $m=15$ | 0.784 | 0.244 | 0.108 | 0.972 | 0.528 | 0 | 0 | 0.78 |
| $m=20$ | 0.808 | 0.432 | 0.084 | 0.972 | 0.628 | 0.028 | 0 | 0.808 |
| $m=25$ | 0.804 | 0.428 | 0.136 | 0.992 | 0.66 | 0.06 | 0 | 0.796 |
| $m=30$ | 0.82 | 0.484 | 0.152 | 0.984 | 0.72 | 0.14 | 0 | 0.796 |
| $m=35$ | 0.828 | 0.492 | 0.18 | 0.988 | 0.74 | 0.156 | 0 | 0.832 |
| $m=39$ | 0.848 | 0.62 | 0.196 | 0.992 | 0.764 | 0.164 | 0 | 0.828 |

The first and second empirical principal component functions each have four extrema, at $-1.5$, 0.5, 2.5 and 10.5, and 0.5, 2.5, 4.5 and 12.5, respectively. Table 1 summarizes the values of $\hat{\pi}$ obtained using 500 bootstrap iterations, and bootstrap sample size $m$ varying among $10, 15, 20, 25, 30, 35, 39$. The value of $\varepsilon$ equals half the length of the interval $\mathcal{J}$.

The tabulated results suggest that the points $-0.5$ and 11.5, and 1.5 and 13.5, correspond to genuine extrema for the first and second empirical principal component functions, respectively. The other extrema are spurious.

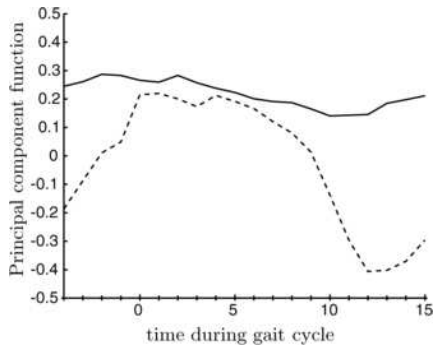

FIG. 1. *Graphs of the first two empirical principal component functions for the hip-movement portion of the child gait data. Solid and dashed lines show the first and second empirical principal component functions, respectively.*

In each instance the largest value of $\hat{\pi}$ among all of the spurious extrema is less than the smallest value of $\hat{\pi}$ among all the genuine extrema.

This information is helpful in choosing the level of smoothing when graphically presenting the principal component curves—the level of smoothing should be sufficiently great to remove the spurious extrema. That is readily achieved, and produces eigenfunction and covariance estimates that reflect the knowledge acquired above.

4.2. *Application to simulated data.* For most of the examples discussed below, synthetic data, representing random functions on the interval $\mathcal{I} = [0,1]$, were generated as follows. The first principal component function was taken to belong to the family

$$(4.1) \quad \psi_1(x) = \psi_{1\lambda}(x) = c_\lambda(2\lambda x^3 - 3\lambda x^2 + 1.5x), \qquad -\infty < \lambda < \infty,$$

where $c_\lambda > 0$ was chosen to ensure that $\int_\mathcal{I} \psi_1^2 = 1$. If $\lambda < 1$ then $\psi_{1\lambda}$ has no extrema or shoulder points on $\mathcal{I}$; if $\lambda = 1$, $\psi_{1\lambda}$ has no extrema and just one shoulder point; if $\lambda > 1$, $\psi_{1\lambda}$ has just two extrema, at the points $[\lambda \pm \{\lambda(\lambda - 1)\}^{1/2}]/(2\lambda)$, and no shoulder point.

This "evolution" of $\psi_{1\lambda}$ takes place smoothly as $\lambda$ is increased, as indicated by the graphs in panel (a) of Figure 2. The other principal component functions, $\psi_k = \psi_{k\lambda}$ for $k \geq 2$, were constructed by orthonormalization from the functions $\{\psi_1(x) : \sin[2\pi(k-1)x], k \geq 2\}$. In our numerical work we used a discrete orthonormalization procedure based on 250 equally-spaced points in $\mathcal{I}$. Then, to build the random function $X$, we multiplied each $\psi_{k\lambda}$ by the square roots of the eigenvalues $\theta_k = k^{-2}$, for $1 \leq k \leq 5$ (we set $\theta_k = 0$ for $k \geq 6$), and also by the $k$th member of an uncorrelated Gaussian random sequence $\eta_1, \eta_2, \ldots$ with zero means and unit variances. In this way we constructed $X = \sum_k \xi_k \psi_{k\lambda}$ [cf. (2.8)], where $\xi_k = \theta_k^{1/2} \eta_k$. In particular, $X$ is a Gaussian process with zero mean.

For this process, numerical results illustrating Theorem 2.1 are readily obtained. They show that the shapes of estimated principal component functions converge to those of true principal component functions, as measured by criteria such as number of extrema and the locations of those points. However, for the sake of brevity we shall not give details of those results.

Panel (b) of Figure 2 graphs, for one value of $m$, the probability density of $\hat{\pi}$ for the principal component functions $\psi_1$ shown in panel (a), respectively. In the case of the solid lines in the figure, the function $\psi_1$ has no extrema and, as predicted by Theorem 3.1, the distribution of $\hat{\pi}$ there is concentrated at relatively low values, that is, toward zero. Considering the dotted lines of the figure, the function $\psi_1$ can be deduced to have two extrema, one a minimum and one a maximum, and, again as suggested by Theorem 3.1, the distribution of $\hat{\pi}$ is concentrated toward the upper end of the unit interval.

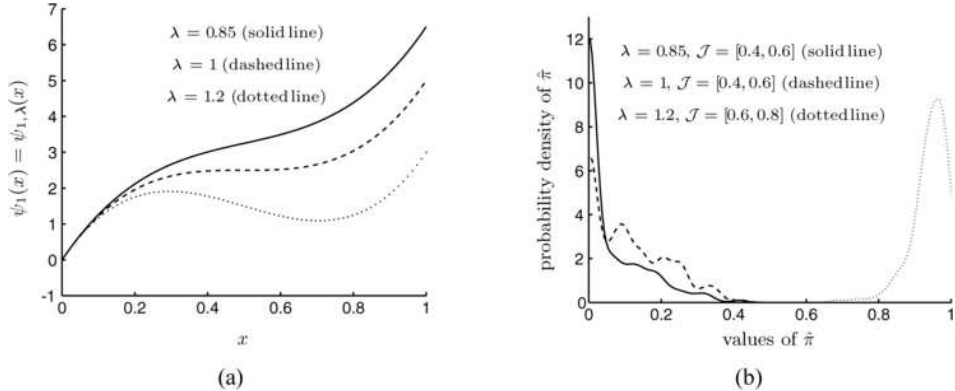

FIG. 2. *Illustration of Theorem* 3.1(i). *In panel* (a) *the function* $\psi_1(x) = \psi_{1\lambda}(x)$, *defined at* (4.1), *is graphed against $x$ for $\lambda = 0.85$, 1 and 1.2. Panel* (b) *shows plots of the probability density of $\hat\pi$, the latter defined at* (3.1), *computed for the first empirical principal component function in the cases of the respective values of $\lambda$. We chose also $\mathcal{J} = [0.4, 0.6]$ (resp. $\mathcal{J} = [0.6, 0.8]$) when $\lambda = 0.85$ or $\lambda = 1$ (resp., $\lambda = 1.2$). The process $X$ was Gaussian, and was constructed as described in the first two paragraphs of this section. The plots in panel* (b) *are for $n = 300$, when $m$ is the integer part of $2n^{0.6}$.*

The dashed lines in Figure 2 show an intermediate case, where $\psi_1$ has a shoulder point and the distribution of $\hat\pi$ is concentrated more toward the centre of the unit interval. Here, and in all the numerical work in this section, sample size was $n = 300$, all results represent averages over 200 synthetic samples of that size, and we used 300 bootstrap iterations for computing $\hat\pi$.

The graphs of the density of $\hat\pi$ in Figure 2 were constructed using kernel methods, in which the kernel was taken to be the standard normal density, and the bandwidth was chosen equal to 0.02, a value suggested by cross-validation. In panel (b) the subinterval $\mathcal{J}$ was centred at 0.5 for the cases $\lambda = 0.85$ and 1. When $\lambda = 1.2$ it was centred at 0.70, close to the local minimum of $\psi$ in $\mathcal{I}$. In each instance, $\mathcal{J}$ was of length 0.2.

Results very similar to those in Figure 2 were obtained in the non-Gaussian case where the variables $\eta_k$ used to construct $X$, four paragraphs above, were taken to have the distribution of a centred value of the absolute value of a standard normal random variable, instead of being standard normal themselves. The probability densities of $\hat\pi$ in the three cases ($\lambda = 0.85$, 1 and 1.2) were almost as well separated as before.

Next we return to the Gaussian case, but alter the definition of $X$ by taking $\psi_1(x) = \psi(x) = (x - 0.8)^3 (x - 0.2)^3$, all other aspects of the construction remaining the same. The function $\psi$ is illustrated in the left-hand panel of Figure 3; it has one local minimum, at 0.5, and two shoulder points, at 0.2 and 0.8. The right-hand panel of the figure shows graphs of the density of $\hat\pi$. For each of those graphs, $\mathcal{J}$ is of width 0.2. The density of $\hat\pi$ has much of its mass toward 0 when $\mathcal{J}$ straddles a shoulder point, but is concentrated close to 1 when $\mathcal{J}$ is centred at the local minimum, reflecting the claims made in Theorem 3.1.

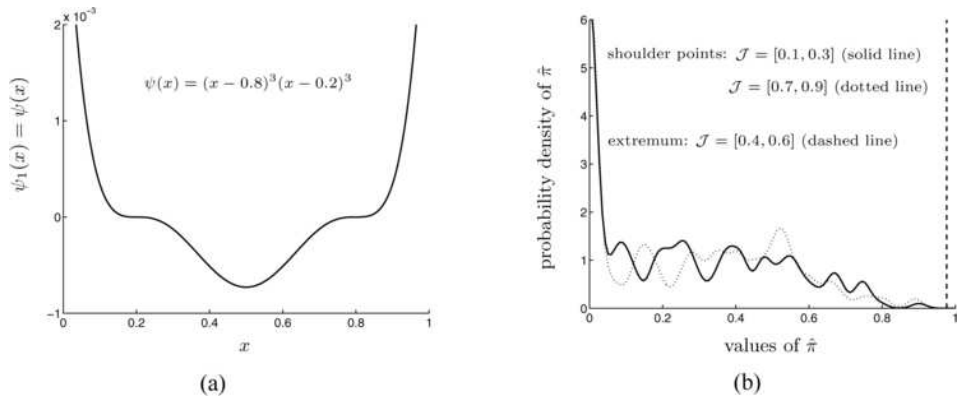

FIG. 3. *Illustration of Theorem* 3.1(ii). *Panel* (a) *graphs the function* $\psi_1(x) = (x - 0.8)^3 (x - 0.2)^3$, *being the first principal component function of a new Gaussian process* $X$. *Panel* (b) *shows three curves, representing the densities of* $\hat{\pi}$, *computed for the first empirical principal component function in different settings. Sample size was* $n = 300$ *and* $m$ *is the integer part of* $2 n^{0.6}$.

Next we construct the process $X$ truncated to 15, rather than 5, nonnull eigenvalues. That is, we take $\lambda_k = k^{-2}$ for $k = 1, \ldots, 15$, and $\lambda_k = 0$ for $k \geq 16$. In this setting, panel (a) of Figure 4 graphs the proportion of variability for each eigenvalue, and suggests that the fifth principal component function is the last one that is likely to be of empirical interest. Panel (b) of Figure 4 graphs the probability density of $\hat{\pi}$ when $\psi_5$ contains either one shoulder point or two extrema. The two curves were obtained using kernel methods with bandwidths 0.02 and 0.01, respectively, and the interval $\mathcal{J}$ was chosen of length 0.05 and centred at 0.5 for the shoulder point, and at 0.7 in the case of the extrema. As expected, the function represented by the solid line is concentrated close to the value 1, but the function represented by the dashed line is supported on a longer interval. This reflects the fact that, in the case of a shoulder point, the empirical evidence for an extremum is very weak.

Next we summarize simulations addressing the local perturbation theory in Section 2.2. We shall use the notation of that section. To address the result (2.10) we constructed the functions $\phi_j$ and $\phi_p$, and the process $Y$, as explained in Section 2.2, with $j = 1$, $p = 2$,

$$\psi_j(x) = 8\sqrt{7}(x - 0.5)^3, \qquad \psi_p(x) = -4\sqrt{5}(x - 0.5)^2.$$

In particular, $\int_{\mathcal{I}} \psi_j^2 = \int_{\mathcal{I}} \psi_p^2 = 1$ and $\int_{\mathcal{I}} \psi_j \psi_p = 0$. The random sequence $\xi_k$, used for the construction of the process $Y$, is defined exactly as described in the second paragraph of this section.

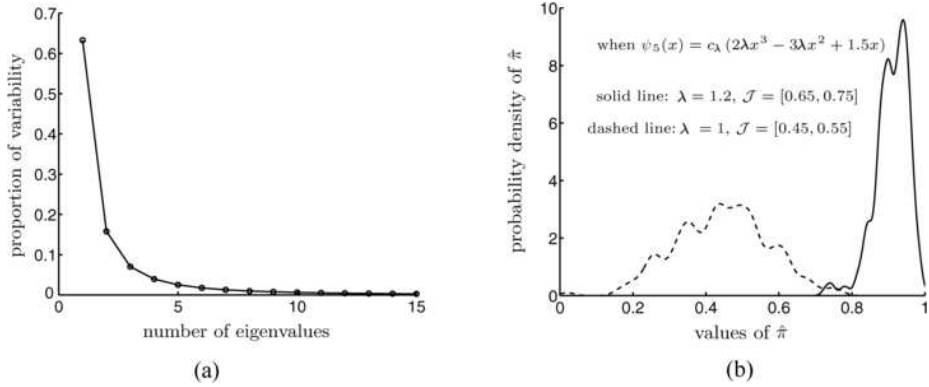

FIG. 4. *Illustration of Theorem* 3.1(iii). *Panel* (a) *graphs the amount of total variation accounted for by each eigenvalue. Panel* (b) *shows two curves, representing the densities of* $\hat{\pi}$, *computed for the fifth empirical principal component function where it possesses a shoulder point (when* $\lambda = 1$), *or an extremum (when* $\lambda = 1.2$). *The interval* $\mathcal{J}$ *was* $[0.45, 0.55]$ *in the case of a shoulder point, and* $[0.65, 0.75]$ *in the case of an extremum.*

Figure 5 illustrates the three cases arising in Theorem 2.2. The results there lend support to the theoretical results in Section 2.2. Indeed, convergence toward 1, or $\frac{1}{2}$, or to a value strictly between 0 and 1, is illustrated by the solid, the dashed and the dashed-dotted curves, respectively. In the last of these three curves the value is close to 0.6 rather than $\frac{1}{2}$.

## 5. Theoretical arguments.

5.1. *Proof of Theorem* 2.1. The following result, which states large-deviation properties in the case of random functions, may be derived using relatively standard arguments. We shall therefore give only an outline proof. Here and below, $C, C_1, C_2, \ldots$ denote constants that may be taken arbitrarily large. We use $1/C, 1/C_1, 1/C_2, \ldots$ to indicate constants whose values can be taken arbitrarily small. Note too that, in results such as (5.1) below, if the result can be proved for any given value of $C$ (say, $C = C_0$), then it also holds for all $C < C_0$. This means that we could equivalently state (5.1) with $1/C$ and $C$, on the left- and right-hand sides, respectively, replaced by $\varepsilon$, and assert that (5.1) holds for all $C, \varepsilon > 0$.

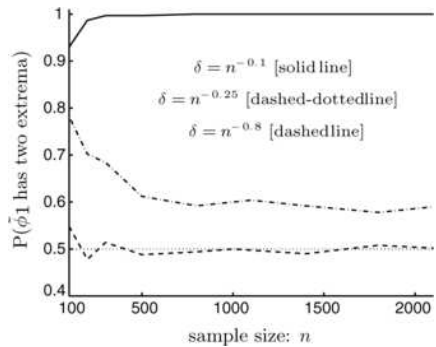

FIG. 5. *Illustration of Theorem* 2.2. *The three curves represent the probabilities that* $\tilde{\phi}_j$ *has two extrema, and no extrema, respectively. The cases* $\delta = n^{-0.1}$, $\delta = n^{-0.25}$ *and* $\delta = n^{-0.8}$ *represent instances where* $n^{1/4}\delta \to \infty$, $n^{1/4}\delta$ *is neither particularly large nor particularly small, and* $n^{1/4}\delta \to 0$, *respectively.*

LEMMA 5.1. *Assume the eigenvalues $\theta_1, \ldots, \theta_{j+1}$ are distinct, and that with probability 1, $X$ has $r \geq 0$ Hölder-continuous derivatives on $\mathcal{I}$, with Hölder exponent $\eta > 0$. If $C > 0$ is given, and if $E|X^{(s)}(a)|^D < \infty$ for $s = 0, \ldots, r$ and $\gamma_r(D, \eta) < \infty$ for a sufficiently large value of $D = D(C) > 0$, then*

$$(5.1) \quad P\left\{\sup_{u \in \mathcal{I}} |\widehat{\psi}_k^{(r)}(u) - \psi_k^{(r)}(u)| > n^{(1/C)-(1/2)} \text{ for } 1 \leq k \leq j\right\} = O(n^{-C}).$$

*Furthermore, if $u_0 \in \mathcal{I}$ and $\eta' \in (0, \eta)$, and if $\gamma_r(D, \eta) < \infty$ for sufficiently large $D = D(C, \eta') > 0$, then*

$$(5.2) \quad P\left\{\sup_{u \in \mathcal{I}} |\widehat{\psi}_k^{(r)}(u) - \widehat{\psi}_k^{(r)}(u_0) - \{\psi_k^{(r)}(u) - \psi_k^{(r)}(u_0)\}| \right.$$
$$\left. > n^{(1/C)-(1/2)}|u - u_0|^{\eta'} \text{ for } 1 \leq k \leq j\right\} = O(n^{-C}).$$

OUTLINE PROOF OF LEMMA 5.1. If $L$ is a function on $\mathcal{I}^2$, put $\|L\|^2 = \int_{\mathcal{I}^2} L^2$, $\|L(u, \cdot)\|^2 = \int_{\mathcal{I}} L(u, v)^2 \, dv$, $\|L\|_{\sup} = \sup_{u \in \mathcal{I}} \|L(u, \cdot)\|$ and $L^{(r)}(u, v) = (\partial/\partial u)^r L(u, v)$. It may be proved that

$$(5.3) \quad \begin{aligned} \hat{\theta}_j \{\widehat{\psi}_j^{(r)}(u) - \psi_j^{(r)}(u)\} \\ = \int_{\mathcal{I}} K^{(r)}(u, v) \{\widehat{\psi}_j(v) - \psi_j(v)\} \, dv \\ + \int_{\mathcal{I}} \{\widehat{K}^{(r)}(u, v) - K^{(r)}(u, v)\} \widehat{\psi}_j(v) \, dv - (\hat{\theta}_j - \theta_j) \psi_j^{(r)}(u), \end{aligned}$$

which entails

$$(5.4) \quad \begin{aligned} \max(0, \theta_j - |\hat{\theta}_j - \theta_j|) |\widehat{\psi}_j^{(r)}(u) - \psi_j^{(r)}(u)| \\ \leq \|K^{(r)}(u, \cdot)\| \|\widehat{\psi}_j - \psi_j\| \\ + \|\widehat{K}^{(r)}(u, \cdot) - K^{(r)}(u, \cdot)\| + |\hat{\theta}_j - \theta_j| |\psi_j^{(r)}(u)|. \end{aligned}$$

Put $\delta_j = \inf_{k \leq j}(\theta_k - \theta_{k+1})$ and $\Delta = \|\widehat{K} - K\|$. It may be proved from results of Bhatia, Davis and McIntosh [2] that $\sup_{j \geq 1} |\hat{\theta}_j - \theta_j| \leq \Delta$, and

that if $\Delta \leq \frac{1}{2}\delta_j$ then $\|\widehat{\psi}_j - \psi_j\| \leq 4\Delta/\delta_j$. On the other hand, if $\Delta > \frac{1}{2}\delta_j$ then $4\Delta/\delta_j > 2 \geq \|\widehat{\psi}_j - \psi_j\|$, since $\|\psi_j\| = \|\widehat{\psi}_j\| = 1$. Therefore, $\|\widehat{\psi}_j - \psi_j\| \leq 4\Delta/\delta_j$ in all cases. When combined with (5.4), these results imply that

$$
\begin{aligned}
\max(0, \theta_j - \Delta)&\|\widehat{\psi}_j^{(r)} - \psi_j^{(r)}\|_\infty \\
&\leq \Delta\{4\|K^{(r)}\|_{\sup}\delta_j^{-1} + \|\psi_j^{(r)}\|_\infty\} + \|\widehat{K}^{(r)} - K^{(r)}\|_{\sup}.
\end{aligned}
\tag{5.5}
$$

The assumptions in Lemma 5.1 imply that $\theta_j, \delta_j > 0$ and $\|K^{(r)}\|_{\sup}$, $\|\psi_j^{(r)}\|_\infty < \infty$. We shall show that if $\gamma_r(D, \eta) < \infty$ for sufficiently large $D = D(C) > 0$, then

$$
P\{\|\widehat{K}^{(r)} - K^{(r)}\|_{\sup} > n^{(1/C)-(1/2)}\} = O(n^{-C}).
\tag{5.6}
$$

A similar but simpler argument shows that the same bound applies if $\|\widehat{K}^{(r)} - K^{(r)}\|_{\sup}$ is replaced by $\Delta$. These results and (5.5) imply the bound (5.1) in the case $k = j$. Other values of $k$ may be treated similarly.

Note that $\widehat{K} = \Delta_1 - \Delta_2$, where

$$\Delta_1(u,v) = \frac{1}{n}\sum_{i=1}^n \{X_i(u) - \mu(u)\}\{X_i(v) - \mu(v)\},$$

$$\Delta_2(u,v) = \{\bar{X}(u) - \mu(u)\}\{\bar{X}(v) - \mu(v)\}.$$

Therefore it suffices to prove the version of (5.6) that arises if we replace $\widehat{K}^{(r)} - K^{(r)}$ there by either $\Delta_3^{(r)} \equiv \Delta_1^{(r)} - K^{(r)}$ or $\Delta_2^{(r)}$. We shall treat only the first of these cases; the second is simpler. That is, we shall prove that

$$
P\left\{\sup_{u \in \mathcal{I}} \|\Delta_3^{(r)}(u, \cdot)\| > n^{(1/C)-(1/2)}\right\} = O(n^{-C}).
\tag{5.7}
$$

First we derive the version of this result when the supremum is taken only over $u \in \mathcal{I}(C_1)$, denoting the set of $n^{C_1}$ equally-spaced points in $\mathcal{I}$, where $C_1 > 0$ can be arbitrarily large. See (5.9) below. Noting that $\Delta_3^{(r)}(u,v)$ equals the mean of a sum of independent and identically distributed random variables, and using Rosenthal's and Markov's inequalities, we may prove that if $C_2, C_3 > 0$ are given, and if

$$
\sup_{u \in \mathcal{I}} E\{|X^{(r)}(u)|^D\} < \infty,
\tag{5.8}
$$

then

$$\sup_{u \in \mathcal{I}} E\{\|n^{1/2}\Delta_3^{(r)}(u, \cdot)\|^D\} < \infty.$$

Therefore, by Markov's inequality,

$$\sup_{u \in \mathcal{I}} P\{\|\Delta_3^{(r)}(u, \cdot)\| > n^{(1/C_2)-(1/2)}\} = O(n^{-D/C_2}).$$

It follows that if $C_2, C_3 > 0$ are given, and if $D = D(C_2, C_3) > 0$ is sufficiently large, then

$$\sup_{u \in \mathcal{I}} P\{\|\Delta_3^{(r)}(u, \cdot)\| > n^{(1/C_2)-(1/2)}\} = O(n^{-C_3}).$$

Hence, if $C_1 > 0$ is any other constant, and if (5.8) holds for sufficiently large $D = D(C_2, C_3, C_1) > 0$, then

$$(5.9) \qquad P\left\{\sup_{u \in \mathcal{I}(C_1)} \|\Delta_3^{(r)}(u, \cdot)\| > n^{(1/C_2)-(1/2)}\right\} = O(n^{-C_3}).$$

If $C_4 > 0$ is given, if $C_1 > 0$ is sufficiently large, and if $\gamma_r(D, \eta) < \infty$ for sufficiently large $D$, then

$$E\left\{\sup_{u \in \mathcal{I}} \|\Delta_3^{(r)}(u, \cdot) - \Delta_3^{(r)}(u', \cdot)\|\right\} = O(n^{-C_4}),$$

where $u'$ denotes the nearest point in $\mathcal{I}(C_1)$ to $u \in \mathcal{I}$. From this bound and Markov's inequality, it follows that if $C_2, C_3 > 0$ are given, and $C_1 = C_1(C_2, C_3) > 0$ is chosen sufficiently large,

$$(5.10) \quad P\left\{\sup_{u \in \mathcal{I}} \|\Delta_3^{(r)}(u, \cdot) - \Delta_3^{(r)}(u', \cdot)\| > n^{(1/C_2)-(1/2)}\right\} = O(n^{-C_3}).$$

The desired result (5.7) follows from (5.9) and (5.10), on taking $C_2 = C_3 > C$.

To derive (5.2) we start from (5.3), which implies that

$$\hat{\theta}_j[\widehat{\psi}_j^{(r)}(u) - \widehat{\psi}_j^{(r)}(u_0) - \{\psi_j^{(r)}(u) - \psi_j^{(r)}(u_0)\}]$$

$$= \int_{\mathcal{I}} \{K^{(r)}(u, v) - K^{(r)}(u_0, v)\}\{\widehat{\psi}_j(v) - \psi_j(v)\} \, dv$$

$$(5.11) \qquad + \int_{\mathcal{I}} [\widehat{K}^{(r)}(u, v) - K^{(r)}(u, v)$$

$$\qquad - \{\widehat{K}^{(r)}(u_0, v) - K^{(r)}(u_0, v)\}]\widehat{\psi}_j(v) \, dv$$

$$- (\hat{\theta}_j - \theta_j)\{\psi_j^{(r)}(u) - \psi_j^{(r)}(u_0)\}.$$

The conditions of the lemma imply that $\psi_j^{(r)}$ is Hölder-continuous with any exponent $\eta' < \eta$. This property, (5.11) and the arguments leading to Lemma 5.1 then give (5.2). $\square$

PROOF OF THEOREM 2.1. For given $\eta > 0$, statements made below hold provided $\gamma_2(D, \eta) < \infty$ for sufficiently large $D > 0$. We shall, however, omit that qualification.

Since the function $|\psi_j'|$ is bounded away from zero uniformly on each of the intervals mentioned in part (a) of the theorem, then part (a) follows directly from (5.1). Let $\mathcal{M}_{k\ell}$ denote the $\varepsilon$-neighborhood of $u_{k\ell}$, mentioned in parts (b) and (c) of the theorem. Using first (5.2) and then (2.4) it may be proved that, for $r = 1, 2$,

$$
\begin{aligned}
\widehat{\psi}_k^{(r)}(u) &= \widehat{\psi}_k^{(r)}(u_{k\ell}) + \psi_k^{(r)}(u) - \psi_k^{(r)}(u_{k\ell}) + O_p(n^{(1/C)-(1/2)}|u - u_{k\ell}|^{\eta'}) \\
&= \widehat{\psi}_k^{(r)}(u_{k\ell}) - \psi_k^{(r)}(u_{k\ell}) \\
&\quad + (\partial/\partial u)^r A_{k\ell}(u - u_{k\ell})^{r_{k\ell}} + O(|u - u_{k\ell}|^{r_{k\ell}+\eta-r}) \\
&\quad + O_p(n^{(1/C)-(1/2)}|u - u_{k\ell}|^{\eta'}),
\end{aligned}
$$
(5.12)

uniformly in $u \in \mathcal{M}_{k\ell}$. Therefore, if $\widehat{\mathcal{S}}_{k\ell}$ denotes the set of all solutions, $u = \hat{u}_{k\ell}$, of the equation $\widehat{\psi}_k'(u) = 0$ in $\mathcal{M}_{k\ell}$, then, taking $r = 1$ in (5.12) and noting that $\psi_k'(u_{k\ell}) = 0$, we deduce first that for each $\delta_1 > 0$, with probability converging to 1 as $n \to \infty$, all elements of $\widehat{\mathcal{S}}_{k\ell}$ lie within $n^{-[1/\{2(r_{k\ell}-1)\}]+\delta_1}$ of $u_{k\ell}$; and then, that for some $\delta_2 > 0$,

$$
(5.13) \quad \hat{u}_{k\ell} - u_{k\ell} = \{-(A_{k\ell} r_{k\ell})^{-1} \widehat{\psi}_k'(u_{k\ell})\}^{1/(r_{k\ell}-1)} + O_p(n^{-[1/\{2(r_{k\ell}-1)\}]-\delta_2}),
$$

uniformly in $\hat{u}_{k\ell} \in \widehat{\mathcal{S}}_{k\ell}$.

Substituting the expansion at (5.13) into the version of (5.12) for $r = 2$ and $u = \hat{u}_{k\ell}$, and writing $\widehat{\phi}_{k\ell}$ for $\widehat{\psi}_k^{(2)}(u_{k\ell}) - \psi_k^{(2)}(u_{k\ell})$, we deduce that for some $\delta_3 > 0$,

$$
\begin{aligned}
\widehat{\psi}_k^{(2)}(\hat{u}_{k\ell}) &= \widehat{\phi}_{k\ell} + r_{k\ell}(r_{k\ell}-1) A_{k\ell} (\hat{u}_{k\ell} - u_{k\ell})^{r_{k\ell}-2} \\
&\quad + O_p(|\hat{u}_{k\ell} - u_{k\ell}|^{r_{k\ell}+\eta-2} + n^{(1/C)-(1/2)}|\hat{u}_{k\ell} - u_{k\ell}|^{\eta'}) \\
&= \widehat{\phi}_{k\ell} + r_{k\ell}(r_{k\ell}-1) A_{k\ell} \{-(A_{k\ell} r_{k\ell})^{-1} \widehat{\psi}_k'(u_{k\ell})\}^{(r_{k\ell}-2)/(r_{k\ell}-1)} \\
&\quad + O_p(n^{-[(r_{k\ell}-2)/\{2(r_{k\ell}-1)\}]-\delta_3}),
\end{aligned}
$$
(5.14)

uniformly in $\hat{u}_{k\ell} \in \widehat{\mathcal{S}}_{k\ell}$, where (here and below) we interpret $x^{(r_{k\ell}-2)/(r_{k\ell}-1)}$ as the $(r_{k\ell}-2)$nd power of $x^{1/(r_{k\ell}-1)}$; it is undefined if $r_{k\ell}$ is odd and $x < 0$. Of course, when $r_{k\ell} = 2$ we interpret $x^{(r_{k\ell}-2)/(r_{k\ell}-1)}$ as 1.

Next we derive the limiting distribution of $(\widehat{\psi}_k'(u_{k\ell}), \widehat{\phi}_{k\ell})$. Observe from (5.3) that

$$
\begin{aligned}
\hat{\theta}_k \{\widehat{\psi}_k^{(r)}(u_{k\ell}) - \psi_k^{(r)}(u_{k\ell})\} \\
= \int_{\mathcal{I}} K^{(r)}(u_{k\ell}, v) \{\widehat{\psi}_k(v) - \psi_k(v)\} dv \\
+ \int_{\mathcal{I}} \{\widehat{K}^{(r)}(u_{k\ell}, v) - K^{(r)}(u_{k\ell}, v)\} \widehat{\psi}_k(v) dv - (\hat{\theta}_k - \theta_k) \psi_k^{(r)}(u_{k\ell}).
\end{aligned}
$$

(We shall need only $r = 1, 2$.) Using arguments of Dauxois, Pousse and Romain [12], it may be proved that

$$\hat{\theta}_k - \theta_k = U_k + o_p(n^{-1/2}),$$

$$\int_{\mathcal{I}} \{\widehat{K}^{(r)}(u_{k\ell}, v) - K^{(r)}(u_{k\ell}, v)\}\widehat{\psi}_k(v)\, dv = V_{k\ell}^{(r)} + o_p(n^{-1/2}),$$

$$\int_{\mathcal{I}} K^{(r)}(u_{k\ell}, v)\{\widehat{\psi}_k(v) - \psi_k(v)\}\, dv = \sum_{p=1}^{\nu} \theta_p W_{kp} \psi_p^{(r)}(u_{k\ell}) + n^{-1/2} R_{k\ell 1}(r, \nu),$$

where

$$U_k = \iint_{\mathcal{I}^2} \{\widehat{K}(u, v) - K(u, v)\}\psi_k(u)\psi_k(v)\, du\, dv,$$

$$V_{k\ell}^{(r)} = \int_{\mathcal{I}} \{\widehat{K}^{(r)}(u_{k\ell}, v) - K^{(r)}(u_{k\ell}, v)\}\psi_k(v)\, dv,$$

$$W_{kp} = (\theta_k - \theta_p)^{-1} \iint_{\mathcal{I}^2} \{\widehat{K}(u, v) - K(u, v)\}\psi_k(u)\psi_p(v)\, du\, dv$$

and, for each $\zeta > 0$ and for $m = 1$,

(5.15) $$\lim_{\nu \to \infty} \limsup_{n \to \infty} P\{|R_{k\ell m}(r, \nu)| > \zeta\} = 0.$$

The results in the previous paragraph imply that

$$\theta_k \{\widehat{\psi}_k^{(r)}(u_{k\ell}) - \psi_k^{(r)}(u_{k\ell})\}$$

$$= \sum_{p=1}^{\nu} \theta_p W_{kp} \psi_p^{(r)}(u_{k\ell}) + V_{k\ell}^{(r)} - U_k \psi_{k\ell}^{(r)}(u_{k\ell}) + n^{-1/2} R_{k\ell 2}(r, \nu)$$

$$= \frac{1}{n} \sum_{i=1}^{n} B_i(r, \nu) + n^{-1/2} R_{k\ell 3}(r, \nu),$$

where

$$B_i(r, \nu) = \sum_{p=1}^{\nu} \theta_p \psi_p^{(r)}(u_{k\ell})(\theta_k - \theta_p)^{-1}$$

$$\times \iint_{\mathcal{I}^2} \{X_i(u)X_i(v) - EX_i(u)X_i(v)\}\psi_k(u)\psi_p(v)\, du\, dv$$

$$+ \int_{\mathcal{I}} \{X_i^{(r)}(u_{k\ell})X_i(v) - EX_i^{(r)}(u_{k\ell})X_i(v)\}\psi_k(v)\, dv$$

$$- \psi_k^{(r)}(u_{k\ell}) \iint_{\mathcal{I}^2} \{X_i(u)X_i(v) - EX_i(u)X_i(v)\}\psi_k(u)\psi_k(v)\, du\, dv$$

$$+ n^{-1/2} R_{k\ell}(r, \nu)$$

and $R_{k\ell 2}(r,\nu)$ and $R_{k\ell 2}(r,\nu)$ satisfy (5.15). [In the definition of $B_i(r,\nu)$ we have, in a slight abuse of notation, used $X_i$ for $X_i - E(X_i)$; there is no loss of generality in assuming that $E(X_i) = 0$.] Taking $r = 1$, in which case $\psi_k^{(r)}(u_{k\ell}) = 0$ and $r = 2$, for which $\widehat{\phi}_{k\ell} = \widehat{\psi}_k^{(r)}(u_{k\ell}) - \psi_k^{(r)}(u_{k\ell})$, we deduce that

$$(5.16) \quad \theta_k(\widehat{\psi}_k'(u_{k\ell}), \widehat{\phi}_{k\ell}) = \frac{1}{n}\sum_{i=1}^n (B_i(1,\nu), B_i(1,\nu)) + n^{-1/2} R_{k\ell}(\nu),$$

where

$$(5.17) \quad \lim_{\nu\to\infty}\limsup_{n\to\infty} P\{|R_{k\ell}(\nu)| > \zeta\} = 0,$$

and on this occasion $|\cdot|$ denotes the Euclidean norm.

The random two-vectors $(B_i(1,\nu), B_i(2,\nu))$ are independent and identically distributed with zero means and finite variances. Therefore conventional arguments show that $n^{-1/2}\sum_i (B_i(1,\nu), B_i(2,\nu))$, and hence also $n^{1/2}(\widehat{\psi}_k'(u_{k\ell}), \widehat{\phi}_{k\ell})$ has an asymptotic bivariate normal distribution with finite variance.

Return to (5.14); multiply throughout by $n^{(r_{k\ell}-2)/\{2(r_{k\ell}-1)\}}$; substitute

$$\frac{1}{n\theta_k}\sum_{i=1}^n (B_i(1,\nu), B_i(2,\nu))$$

for $(\widehat{\psi}_k'(u_{k\ell}), \widehat{\phi}_{k\ell})$ when a factor of the latter appears on the right-hand side of (5.14); and note (5.16) and (5.17), obtaining

$$n^{(r_{k\ell}-2)/\{2(r_{k\ell}-1)\}}\widehat{\psi}_k^{(2)}(\hat{u}_{k\ell})$$

$$= n^{-r_{k\ell}/\{2(r_{k\ell}-1)\}}\theta_k^{-1}\sum_{i=1}^n B_i(2,\nu)$$

$$+ r_{k\ell}(r_{k\ell}-1)A_{k\ell}\left\{-(A_{k\ell}r_{k\ell})^{-1}\frac{1}{n^{1/2}}\sum_{i=1}^n B_i(1,\nu)\right\}^{(r_{k\ell}-2)/(r_{k\ell}-1)}$$

$$+ o_p(1),$$

uniformly in $\hat{u}_{k\ell} \in \widehat{\mathcal{S}}_{k\ell}$. The series $\sum_i B_i(2,\nu)$ equals $O_p(n^{1/2})$, which, when multiplied by $n^{-r_{k\ell}/\{2(r_{k\ell}-1)\}}$, equals $o_p(1)$. Therefore,

$$n^{(r_{k\ell}-2)/\{2(r_{k\ell}-1)\}}\widehat{\psi}_k^{(2)}(\hat{u}_{k\ell})$$

$$= r_{k\ell}(r_{k\ell}-1)A_{k\ell}\left\{-(A_{k\ell}r_{k\ell})^{-1}\frac{1}{n^{1/2}}\sum_{i=1}^n B_i(1,\nu)\right\}^{(r_{k\ell}-2)/(r_{k\ell}-1)} + o_p(1)$$

$$= r_{k\ell}(r_{k\ell}-1)A_{k\ell}N_{k\ell}(n)^{(r_{k\ell}-2)/(r_{k\ell}-1)} + o_p(1),$$

uniformly in $\hat{u}_{k\ell} \in \widehat{\mathcal{S}}_{k\ell}$, where the random variable $N_{k\ell}(n)$ converges in distribution to a random variable $N_{k\ell}$ with the $N(0, \sigma_{k\ell}^2)$ distribution.

When $r_{k\ell} \geq 2$ is even, the sign of $r_{k\ell}(r_{k\ell} - 1)A_{k\ell}N_{k\ell}^{(r_{k\ell}-2)/(r_{k\ell}-1)}$ is always the same as the sign of $A_{k\ell}$. [Recall the interpretation, just below (5.14), of $x^{(r_{k\ell}-2)/(r_{k\ell}-1)}$.] Therefore, when $r_{k\ell} \geq 2$, all the solutions of $\widehat{\psi}_k' = 0$ that lie in the vicinity of $u_{k\ell}$ give, with probability converging to 1 as $n \to \infty$, local maxima if $A_{k\ell} < 0$ or local minima if $A_{k\ell} > 0$. It follows that, if $r_{k\ell} \geq 2$ is even, with probability converging to 1 as $n \to \infty$, $\widehat{\psi}_k$ has just one extremum in the neighborhood of $u_{k\ell}$, with the same parity as the extremum of $\psi_k$ at $u_{k\ell}$. This proves part (b) of the theorem.

On the other hand, when $r_{k\ell} \geq 3$ is odd, the random variable $N_{k\ell}^{1/(r_{k\ell}-1)}$ is well-defined if and only if $N_{k\ell} > 0$, in which case it takes either of two values, that is plus or minus the $(r_{k\ell} - 1)$st root of $N_{k\ell}$. Correspondingly, the probability that $\widehat{\mathcal{S}}_{k\ell}$ is empty converges to $P(N_{k\ell} < 0) = \frac{1}{2}$, and the probability that $\widehat{\mathcal{S}}_{k\ell}$ contains just two elements converges to $P(N_{k\ell} > 0)$. The signs of $N_{k\ell}^{(r_{k\ell}-2)/(r_{k\ell}-1)}$ that correspond to the two values of $N_{k\ell}^{1/(r_{k\ell}-1)}$ are positive and negative, respectively. Correspondingly, the signs of $N_{k\ell}(n)^{(r_{k\ell}-2)/(r_{k\ell}-1)}$ are positive and negative, and so when $\widehat{\mathcal{S}}_{k\ell}$ contains two elements, one of them must give a local maximum and the other a local minimum. This proves parts (c) and (d) of the theorem.

When $r_{k\ell}$ is even, that is when $u_{k\ell}$ gives a local extremum of $\psi_k$, the first term on the right-hand side of (5.13) is always well defined, since $x^{1/(r_{k\ell}-1)}$ is well defined for all real $x$. In particular, part (e) of Theorem 2.1 follows directly from (5.13), since $n^{1/2}\widehat{\psi}_k'(u_{k\ell})$ is asymptotically normally distributed.

When $r_{k\ell}$ is odd, that is when $u_{k\ell}$ gives a point of inflection of $\psi_k$, the first term on the right-hand side of (5.13) is well defined only when $-(A_{k\ell}r_{k\ell})^{-1}\widehat{\psi}_k'(u_{k\ell}) \geq 0$. The results discussed two paragraphs above now imply that, conditional on $\widehat{\psi}$ having two extrema given generically by $\hat{u}_{k\ell}$, they satisfy

$$n^{1/\{2(r_{k\ell}-1)\}}(\hat{u}_{k\ell} - u_{k\ell}) = N_{k\ell}(n)^{1/(r_{k\ell}-1)} + o_p(1),$$

where $N_{k\ell}$ is conditioned to be positive and the two solutions correspond to the positive and negative roots on the right-hand side. This gives part (f) of Theorem 2.1. □

5.2. *Proof of Theorem* 3.1. The proof closely parallels that of Theorem 2.1, being based on the following bootstrap version of Lemma 5.1.

LEMMA 5.2. *Assume the conditions of Lemma* 5.1. *If* $C > 0$ *is given, and if* $\gamma_r(D, \eta) < \infty$ *for a sufficiently large value of* $D = D(C) > 0$, *then there* exists a set $\mathcal{R}$ of realisations of $\mathcal{X}$, with probability $P(\mathcal{R}) = 1 - O(n^{-C})$, such that whenever $\mathcal{X} \in \mathcal{R}$,

$$P\left\{\sup_{u \in \mathcal{I}} |(\widehat{\psi}_k^*)^{(r)}(u) - \widehat{\psi}_k^{(r)}(u)| > m^{(1/C) - (1/2)} \text{ for } 1 \leq k \leq j \Big| \mathcal{X}\right\} \leq \text{const.} m^{-C},$$

where, here and below, "const." is nonrandom. Furthermore, if $u_0 \in \mathcal{I}$ and $\eta' \in (0, \eta)$, and if $\gamma_r(D, \eta) < \infty$ for sufficiently large $D = D(C, \eta') > 0$, then whenever $\mathcal{X} \in \mathcal{R}$,

$$P\bigg\{\sup_{u \in \mathcal{I}} |(\widehat{\psi}_k^*)^{(r)}(u) - (\widehat{\psi}_k^*)^{(r)}(u_0) - \{\widehat{\psi}_j^{(r)}(u) - \widehat{\psi}_j^{(r)}(u_0)\}|$$
$$> m^{(1/C)-(1/2)} |u - u_0|^{\eta'} \text{ for } 1 \leq k \leq j \Big| \mathcal{X}\bigg\} \leq \text{const.} m^{-C}.$$

Let $\mathcal{M}_{k\ell}$ be as in the proof of Theorem 2.1, and write $\widehat{\mathcal{S}}_{k\ell}^*$ for the set of all solutions, $u = \hat{u}_{k\ell}^*$, of the equation $(\widehat{\psi}_k^*)'(u) = 0$ in $\mathcal{M}_{k\ell}$. Parallelling arguments in Section 5.1 we may prove the following analogues of (5.13) and (5.14): for some $\delta > 0$,

(5.18)
$$\hat{u}_{k\ell}^* - u_{k\ell} = [(A_{k\ell} r_{k\ell})^{-1} \{(\widehat{\psi}_k^*)'(u_{k\ell}) - \psi_k'(u_{k\ell})\}]^{1/(r_{k\ell}-1)}$$
$$+ m^{-[1/\{2(r_{k\ell}-1)\}]} R_1^*(\hat{u}_{k\ell}^*),$$

(5.19)
$$(\widehat{\psi}_k^*)^{(2)}(\hat{u}_{k\ell}^*) = r_{k\ell}(r_{k\ell} - 1) A_{k\ell}$$
$$\times [(A_{k\ell} r_{k\ell})^{-1} \{(\widehat{\psi}_k^*)'(u_{k\ell}) - \psi_k'(u_{k\ell})\}]^{(r_{k\ell}-2)/(r_{k\ell}-1)}$$
$$+ m^{-[(r_{k\ell}-2)/\{2(r_{k\ell}-1)\}]} R_2^*(\hat{u}_{k\ell}^*),$$

where, on a set $\mathcal{R}(C_1)$ of realizations of $\mathcal{X}$ that satisfies $P(\mathcal{R}) = 1 - O(n^{-C})$,

$$P\left\{\sup_{\hat{u}_{k\ell}^* \in \widehat{\mathcal{S}}_{k\ell}^*} |R_2^*(\hat{u}_{k\ell}^*)| > n^{-\delta} \Big| \mathcal{X}\right\} \leq n^{-C}$$

for all $\mathcal{X} \in \mathcal{R}$.

Write

$$(\widehat{\psi}_k^*)'(u_{k\ell}) - \psi_k'(u_{k\ell}) = (\widehat{\psi}_k^*)'(u_{k\ell}) - \widehat{\psi}_k'(u_{k\ell}) + \widehat{\psi}_k'(u_{k\ell}) - \psi_k'(u_{k\ell})$$

in each of (5.18) and (5.19), and note that the variables

$$Z = n^{1/2} \{\widehat{\psi}_k'(u_{k\ell}) - \psi_k'(u_{k\ell})\}, \qquad Z^* = m^{1/2} \{(\widehat{\psi}_k^*)'(u_{k\ell}) - \widehat{\psi}_k'(u_{k\ell})\}$$

are each asymptotically normal distributed as $N(0, \tau^2)$ say, with $Z^*$ having this weak limit conditional on $\mathcal{X}$ (and hence also conditional on $Z$). Theorem 3.1 follows from these properties, on doing little more than retracing the arguments leading to Theorem 2.1.

For example, to obtain the last part of (c)$'$ in the case where $r_{k\ell}$ is odd, use Kolmogorov's extension theorem to write $Z = Z_1 + o_p(1)$ and $Z^* = Z_2 + o_p(1)$, where $Z_1$ and $Z_2$ have exactly $\mathrm{N}(0, \tau^2)$ distributions, $Z_1$ is measurable in the sigma-field generated by $\mathcal{X}$, and $Z_2$ is $\mathrm{N}(0, \tau^2)$ conditional on $\mathcal{X}$. Assume, without loss of generality, that $A_{k\ell} > 0$; the case $A_{k\ell} < 0$ may be treated similarly. Then, in view of (5.18), the probability, conditional on $\mathcal{X}$, that $\widehat{\mathcal{S}}_{k\ell}^*$ is nonempty (and contains exactly one element) equals

$$(5.20) \qquad P(\rho^{1/2} Z_1 + Z_2 > 0 | Z_1) + o_p(1),$$

where $\rho = \lim(m/n)$. If $\rho = 0$ then the quantity at (5.20) converges in probability to $\frac{1}{2}$; and if $0 < \rho \leq 1$ then it equals $\Phi(\rho^{1/2} Z_3) + o_p(1)$, where $Z_3 = Z_1/\tau$ and has a $\mathrm{N}(0, 1)$ distribution.

Centre for Mathematics
 and its Applications
Australian National University
Canberra ACT 0200
Australia

MODAL'X
Université Paris X-Nanterre
Bat G
200 avenue de la république
9200 Paris
France
E-mail: cvial@u-paris10.fr